\documentclass[12pt]{amsart}
\usepackage{amssymb}
\usepackage{amsmath}
\usepackage{mathtools, float}
\usepackage{stmaryrd}
\usepackage[english]{babel}
\usepackage[utf8]{inputenc}
\usepackage[T1]{fontenc}
\usepackage{enumerate}
\usepackage{graphicx,pinlabel,xcolor, fullpage}
\usepackage{adjustbox}
\usepackage{array,multirow}
\usepackage{float}
\usepackage{oldgerm}

\theoremstyle{plain}

\newtheorem{theorem*}{Theorem}

\theoremstyle{definition}

\newtheoremstyle{case}{}{}{}{}{}{:}{ }{}
\theoremstyle{case}

\newcommand{\QED}{\hfill $\blacksquare$}

\newcommand{\bdot}{\boldsymbol{\cdot}}

\begin{document}
\title[A note on a Conjecture of Gao and Zhuang for groups of order $27$]{A note on a Conjecture of Gao and Zhuang for groups of order $27$}
\author{Naveen K. Godara and Siddhartha Sarkar}
\address{Department of Mathematics\\
Indian Institute of Science Education and Research Bhopal\\
Bhopal Bypass Road, Bhauri \\
Bhopal 462 066, Madhya Pradesh\\
India}
\email{naveen16@iiserb.ac.in, sidhu@iiserb.ac.in}

\vspace*{1em}

\begin{abstract}
The small Davenport constant ${\mathsf{d}}(G)$ of a finite group $G$ is defined to be the maximal length of a sequence over $G$ which has no non-trivial product-one subsequence. In this paper, we prove that ${\mathsf{d}}(G) = 6$ for the non-abelian group of order $27$ and exponent $3$ and thereby establish a conjecture by Gao and Zhuang for this group.      
\end{abstract}


\maketitle

\bigskip


\bigskip

\section{Introduction}
\label{introsec}

\noindent Let $G$ be a finite group, written multiplicatively, and ${\mathcal{F}}(G)$ denote the free abelian monoid generated by $G$ as its free basis. A finite unordered {\it sequence} of elements in $G$ can be viewed as an element of ${\mathcal{F}}(G)$, written as $S = g_1 \bdot \dotsc \bdot g_{\ell}$, where $g_1, \dotsc, g_{\ell} \in G$, which are called the {\it terms} of the sequence $S$. The non-negative integer ${\ell}$ corresponding to $S$ is called the {\it length} of the sequence $S$ and is denoted by $\ell(S)$. The sequence $S$ with $\ell(S) = 0$ is the {\it empty} or {\it trivial} sequence, denoted by $1_{{\mathcal{F}}(G)} \in {\mathcal{F}}(G)$, which is the identity element of ${\mathcal{F}}(G)$. For any $g \in G$, we denote
\[
{\mathsf{v}}_{g}(S) := \lvert \{ i \in [1,\ell] ~:~ g_i = g \} \rvert \hspace*{.2in} ({\mathrm{\it{multiplicity}}}~{\mathrm{of}}~ g ~{\mathrm{in}}~ S)
\]
where $[1,\ell]$ denotes the set $\{ 1, \dotsc, \ell \}$ of positive integers. An element $T$ in ${\mathcal{F}}(G)$ is called a {\it subsequence} of a sequence $S$, denoted by $T \mid S$, if ${\mathsf{v}}_{g}(T) \leq {\mathsf{v}}_{g}(S)$ for every $g \in G$.

\smallskip

\noindent For two sequences $S = g_1 \bdot \dotsc \bdot g_{\ell_1}$ and $T = h_1 \bdot \dotsc \bdot h_{\ell_2}$, we write $S \bdot T := g_1 \bdot \dotsc \bdot g_{\ell_1} \bdot h_1 \bdot \dotsc \bdot h_{\ell_2}$ in ${\mathcal{F}}(G)$. Since this operation $\bdot$ is associative in ${\mathcal{F}}(G)$, for sequences $T_1, \dotsc, T_m \in {\mathcal{F}}(G)$, we write $T_1 \bdot T_2 \bdot \dotsc \bdot T_m$ instead of putting parenthesis among them. In case $T \mid S$ in ${\mathcal{F}}(G)$, we use the notation $T^{[-1]} \bdot S$ or $S \bdot T^{[-1]}$ to denote the sequence obtained from $S$ by removing the terms of $T$. For an element $g \in G$, the sequence of length $1$ with only term $g$ is denoted by $(g)$. For a sequence $S = g_1 \bdot \dotsc \bdot g_{\ell}$ of length $\ell \geq 1$ over $G$, we define 
\[
\pi(S) := \big\{ g_{\tau(1)} \dotsc g_{\tau(\ell)} \in G ~:~ \tau ~{\mathrm{is~a~permutation~of~}} [1,\ell] \big\},
\] 
and 
\[  
\Pi(S) := \bigcup_{1_{{\mathcal{F}}(G)} \neq T \in {\mathcal{F}}(G),~ T \mid S} \pi(T),
\]
as subsets of $G$. 

\smallskip

\noindent Since we would be interested in non-abelian groups $G$, it would be often necessary to consider the product $g_1 \dotsc g_{\ell} \in G$ corresponding to a sequence $g_1 \bdot \dotsc \bdot g_{\ell} \in {\mathcal{F}}(G)$. For this, one can consider the free non-abelian monoid ${\mathcal{F}}^{\ast}(G)$ generated by $G$ as a basis and the natural map $[\bdot] : {\mathcal{F}}^{\ast}(G) \rightarrow {\mathcal{F}}(G)$ that extends the identity map $G \rightarrow G$ (as sets). For a sequence $S = g_1 \bdot \dotsc \bdot g_{\ell} \in {\mathcal{F}}(G)$, we use the notation $S^{\ast}$ to denote the corresponding ordered sequence, denoted by the same notation $S^{\ast} = g_1 \bdot \dotsc \bdot g_{\ell}$, which represents the associated element in ${\mathcal{F}}^{\ast}(G)$, keeping the order of the terms as they appear. With these considerations, we have $[S^{\ast}] = S$, and we define $\pi(S^{\ast}) := g_1 \dotsc g_{\ell} \in G$. For more details on these notations and terminologies, see \cite{GGR2013}. 

\smallskip

\noindent Let $S$ be a non-trivial sequence over $G$. We say that $S$ is a {\it product-one sequence} if $1 \in \pi(S)$. Next, we say that $S$ has a {\it product-one subsequence} if $1 \in \Pi(S)$. The sequence $S$ is said to have a {\it central product} if $\pi(S^{\ast}) \in Z(G)$, where $Z(G)$ denotes the center of the group $G$. Notice that $\pi(S)$ is contained in the coset $\pi(S^{\ast}) [G,G]$. 

\bigskip

\noindent The small Davenport constant ${\mathsf{d}}(G)$ is defined to be the maximal integer $\ell$ such that there is a sequence $S$ over $G$ of length $\ell$ that has no non-trivial product-one subsequence; i.e., $1 \not\in \Pi(S)$. The large Davenport constant ${\mathsf{D}}(G)$ of $G$ is defined to be the maximal length of a minimal product-one sequence: this is a product-one sequence that cannot be factored into two non-trivial, product-one subsequences. The Erd\"{o}s-Ginzburg-Ziv constant ${\mathsf{s}}(G)$ is defined to be the least positive integer $k$, such that every sequence $S$ over $G$ of length $\ell(S) \geq k$ has a product-one subsequence $T$ of length $\ell(T) = {\mathrm{exp}}(G)$, where ${\mathrm{exp}}(G)$ denotes the least positive integer $m$ such that $g^m = 1$ for every $g \in G$.

\bigskip

\noindent The problem of finding the Davenport constant ${\mathsf{D}}(G)$ for a finite abelian group $G$ was posed by H. Davenport in 1966 at the Mid-western Conference on Group Theory and Number Theory in Ohio, where he showed that if $K$ is a number field, ${\mathcal{O}}_K$ its number ring, and $G_K$ the corresponding class group, then for each non-unit and irreducible element $x \in {\mathcal{O}}_K$, the ideal generated by $x$ in ${\mathcal{O}}_K$ can be written as a product of at most ${\mathsf{D}}(G_K)$ many prime ideals with multiplicity. 

\bigskip

\noindent For a finite group $G$, it follows that ${\mathsf{d}}(G) + 1 \leq {\mathsf{D}}(G) \leq |G|$, and in particular, if $G$ is finite abelian, then the equality ${\mathsf{d}}(G) + 1 = {\mathsf{D}}(G)$ holds (see Lemma 2.4, \cite{GGR2013}). In 1969, Olson proved the following result:

\bigskip

\subsection{Theorem}\label{Ols} (\cite{ols1969}) If $p$ is a prime and $G$ is a finite abelian $p$-group of the form 
\[
G(p; e_1, \dotsc, e_r) := C_{p^{e_1}} \times \dotsc \times C_{p^{e_r}}, \hspace*{.5in} 1 \leq e_1 \leq \dotsc \leq e_r 
\]
with $r \geq 1$, then ${\mathsf{D}}(G) = 1 + \sum_{i=1}^{r} (p^{e_i} - 1)$. 

\bigskip

\noindent So far, there is no closed formula known for ${\mathsf{D}}(G)$ for an arbitrary abelian group. Recently, the following sharp result was proved by Bhowmik and Schlage-Puchta (\cite{BSP2012}): for a finite abelian group $G$
\[
{\mathsf{D}}(G) \leq 
\begin{cases}
{\mathrm{exp}}(G) + {\frac {|G|}{ {\mathrm{exp}}(G) }} - 1 &\mbox{if } {\mathrm{exp}}(G) \geq {\sqrt{|G|}}, \\
2 {\sqrt{|G|}} - 1 &\mbox{if } {\mathrm{exp}}(G) < {\sqrt{|G|}}.
\end{cases}
\]

\bigskip

\noindent The small Davenport constant has some interesting connection to the Gao constant ${\mathsf{E}}(G)$ for a finite group $G$, which is defined as the smallest positive integer $N$ such that every sequence $S$ over $G$ of length $\ell(S) \geq N$ has a subsequence $T$ with $\ell(T) = |G|$ and $1 \in \pi(T)$.

\bigskip

\noindent In 1961, Erd\"{o}s, Ginzburg, and Ziv (\cite{EGZ1961}) proved that ${\mathsf{E}}(G) \leq 2|G| - 1$ for any finite solvable group $G$, which is known as the Erd\"{o}s-Ginzburg-Ziv theorem. Later on, Olson (\cite{OLS1976}) extended this result for any finite group $G$. Since 1996, there has been an attempt to reduce this upper bound for finite non-cyclic solvable groups: Yuster and Peterson (1984, \cite{PEY1984}) showed that ${\mathsf{E}}(G) \leq 2|G| - 2$, and later Yuster (1988, \cite{YUS1988}) improved this bound to ${\mathsf{E}}(G) \leq 2|G| - r$ if $|G| \geq 600 ((r-1)!)^2$. In 1996, Gao (\cite{GAO1996B}) improved this bound to ${\mathsf{E}}(G) \leq {\frac {11}{6}} |G| - 1$, and later on, Gao and Li (2010, \cite{GLI2010}) made a further improvement of this to ${\mathsf{E}}(G) \leq {\frac {7}{4}} |G| - 1$. In this last article, Gao and Li posed the following conjecture:

\bigskip

\noindent {\bf Conjecture 1.} For any finite group $G$, we have ${\mathsf{E}}(G) \leq {\frac {3}{2}}|G|$.

\bigskip 

\noindent In 2005, Gao and Zhuang (\cite{ZGA2005}) also posed another conjecture:

\bigskip

\noindent {\bf Conjecture 2.} ${\mathsf{E}}(G) = {\mathsf{d}}(G) + |G|$ for any finite group $G$.  

\bigskip

\noindent In 2015, Han proved Conjecture 1 for all finite non-cyclic nilpotent groups and showed that the solvable group $C_p \ltimes C_{pn}$ satisfies Conjecture 2 as well. So far, both of these conjectures are known to be affirmative for the following class of groups:

\bigskip

\noindent (1) $G$ is finite abelian (Gao, \cite{GAO1996}), \\
\noindent (2) $G$ is Dihedral group of order $2n$ and Dicyclic group of order $4n$ with $n \geq 2$ (Bass, \cite{BAS2007}), \\
\noindent (3) $G \cong C_p \ltimes C_{pn}$ with $p$ a prime and $n \geq 1$ any integer (Han, \cite{HAN2015}). 

\bigskip

\noindent Recently, Gao, Li, and Qu (\cite{GLQ2021}) have modified the bound in Conjecture 1: ${\mathsf{E}}(G) \leq {\frac {3}{2}}(|G| - 1)$ for groups of odd order and $|G| > 9$. See {\cite{ABR2023}, \cite{DHA2019}, \cite{YYU2023}} for recent progress on the above conjectures. For conjectures of Gao relating Erd\"{o}s-Ginzburg-Ziv constant ${\mathsf{s}}(G)$ with certain other constants, see \cite{ABR2023}.  

\bigskip 

\noindent In this paper, we consider non-abelian groups $G$ of order $p^3$, where $p$ is a prime. In cases where either $p=2$ or $p \geq 3$ and $G$ contain a cyclic subgroup of index $p$, the conjectures are confirmed (see \cite{BAS2007}, \cite{HAN2015}). So, we assume that $p$ is an odd prime and $G \cong H_{p^3} := C_p \ltimes ( C_p \times C_p )$, which has order $p^3$ and exponent $p$. In this article, we derive some general results for $H_{p^3}$ and finally confirm the Conjecture 2 for the group $G = H_{27}$.   

\bigskip

\noindent We start with the following definition:

\subsection{Definition}\label{EGZ-seq} Let $p$ be an odd prime and $G$ be the non-abelian $p$-group of order $p^3$ and exponent $p$. A sequence $T$ over $G$ is called an {\it Erd\"{o}s-Ginzburg-Ziv sequence} ({\it EGZ sequence}) if $T$ has central product, and $T$ contains a subsequence $T^{\prime} = g_1 \bdot \dotsc \bdot g_{p-1} \bdot g_p$ such that $g_p$ does not commute with the product $g_{\alpha} \dotsc g_{\beta}$ for any $1 \leq \alpha \leq \beta \leq p-1$. If this happens, we call the subsequence $T^{\prime}$ a {\it principal part} of $T$.  

\bigskip

\noindent Notice that the hypothesis in the above definition ensures that the terms of the principal part of $T$ are non-central in $G$ and $\ell(T) \geq p+1$. Next, if $T_1 \mid T_2$ are sequences over $G$ that have central product, and if $T_1$ is an EGZ sequence, then $T_2$ is also an EGZ sequence. Next, we prove that:  

\bigskip

\subsection{Theorem}\label{EGZ-type2-length-p2} Let $p$ be an odd prime and $G$ be the non-abelian group of order $p^3$ and exponent $p$. Let $U = g_1 \bdot \dotsc \bdot g_m$ be a sequence over $G$ of length $m \geq p+1$ that has central product. If $U$ contains an EGZ subsequence, then there exists a permutation $\sigma$ of $[1,m]$ such that $g_{\sigma(1)} \dotsc g_{\sigma(m)} = 1$. In particular, every EGZ sequence over $G$ contains a product-one subsequence. 

\bigskip

\noindent As an application of these results, we then prove that (Theorems \ref{Davenport-H27} and \ref{EGZ-H27}):

\bigskip

\subsection{Theorem}\label{EGZ-final} Let $G$ be the non-abelian group of order $27$ and exponent $3$. Then ${\mathsf{d}}(G) = 6$ and ${\mathsf{E}}(G) = 33$. In particular, ${\mathsf{E}}(G) = {\mathsf{d}}(G) + |G|$.

\bigskip

\noindent We leave the following question:

\bigskip

\noindent {\bf Question.} Prove that if $p$ is an arbitrary odd prime and $G$ is the non-abelian group of order $p^3$ and exponent $p$, then ${\mathsf{d}}(G) = 3p-3$ and the equality ${\mathsf{E}}(G) = {\mathsf{d}}(G) + |G|$ holds.    

\bigskip

\section{Results on sequences in $H_{p^3}$ that contain EGZ subsequences}
\label{Heisenberg_arbitrary_p}

\bigskip

\noindent In this section, we are going to prove some basic results that lead to the results on $H_{27}$. Here, we always assume $p$ is an arbitrary odd prime, and for most of this section, $G$ denotes the non-abelian group of order $p^3$ and exponent $p$, which is given by
\[
G = \Big\langle x, y ~:~ x^p, y^p, [y,x,y], [y,x,x] \Big\rangle,
\]
unless otherwise stated in certain contexts.

\bigskip

\noindent Then, $[G, G] = \langle [x,y] \rangle = Z(G)$ has order $p$ and $G$ has the abelianization $G/{[G, G]} \cong C_p \times C_p$. Denoting by $v := [y,x]$, the generator of $Z(G)$, any two elements $x^{i_1} y^{j_1} v^{k_1}$ and $x^{i_2} y^{j_2} v^{k_2}$ with $0 \leq i_1, j_1, k_1, i_2, j_2, k_2 \leq p-1$ are conjugate if and only if either $(i_1, j_1) = (i_2, j_2) = (0, 0), k_1 = k_2$ (which represents the central elements of $G$) or else $(i_1, j_1) = (i_2, j_2) \neq (0, 0)$. Thus, $G$ has precisely $p^2 + p - 1$ conjugacy classes, where the $p^2 - 1$ conjugacy classes containing the non-central elements are given by
\[
\{ x^i y^j v^k ~:~ 0 \leq k \leq p-1 \} \hspace*{.2in} (0 \leq i,j \leq p-1, (i,j) \neq (0,0)).
\] 

\noindent The computations to derive the Davenport constant of $G$ are more connected to the concept of a $z$-class of $G$, which is a natural generalization of a conjugacy class. 

\bigskip

\subsection{Definition} Let $G$ be any group, and let $g,h \in G$. We call $g$ is $z$-equivalent to $h$ in $G$, denoted by $g \sim_z h$, if there exists an element $x \in G$ such that $x^{-1} C_{G}(g) x = C_{G}(h)$.

\bigskip

\noindent The relation $\sim_z$ is an equivalence relation in $G$. In cases where $g$ and $h$ are conjugate in $G$, they are $z$-equivalent. Consequently, the orbits of the equivalence relation $\sim_z$ in $G$ are unions of certain conjugacy classes of $G$. In particular, the centre $Z(G)$ of $G$ forms a $z$-class in $G$. The following Lemma is easy to prove, and we record it for its usefulness. 

\bigskip

\subsection{Lemma} Let $p$ be an odd prime and $G$ be the non-abelian $p$-group of order $p^3$ and exponent $p$. Then:

\smallskip

\noindent (i) For any $g \in G \setminus Z(G)$, the $z$-class containing $g$ is given by 
\[
{\mathcal{K}}[g] := \big\{ g^\lambda u ~:~ 1 \leq \lambda \leq p-1, u \in Z(G) \big\} = C_G(g) \setminus Z(G).
\]
Thus, for any $g,h \in G$, $g$ centralizes $h$ if and only if either one of them is central, or else ${\mathcal{K}}[g] = {\mathcal{K}}[h]$. In particular, there are precisely $p+2$ $z$-classes of $G$.  

\smallskip

\noindent (ii) If ${\overline{g}}$ (resp. ${\overline{Y}}$) denotes the image of $g \in G$ (resp. $Y$ a subset of $G$) under the natural quotient map $: G \rightarrow G/{Z(G)}$, then for any $g \in G \setminus Z(G)$, the set ${\overline{ {\mathcal{K}}[g] }} \cup \{ {\overline{1}} \}$ is a subgroup of order $p$ in ${\overline{G}} \cong C_p \times C_p$. 

\subsection{Definition}\label{notation-z-classes} For any $0 \leq i \leq p+1$, we denote by ${\mathcal{K}}_i$ the $z$-classes of $G$ as follows:
\[
{\mathcal{K}}_i := 
\begin{cases}
Z(G) &\mbox{if } i=0, \\
{\mathcal{K}}[x y^{i-1}] &\mbox{if } 1 \leq i \leq p, \\
{\mathcal{K}}[y] &\mbox{if } i=p+1.
\end{cases}
\]

\subsection{Lemma}\label{BEN-Lemma} (cf. \cite{ben}) Let $T$ be a sequence over a finite cyclic group of order $n$. If $\ell(T) \geq {\frac {n+1}{2}}$ and $T$ contains no non-empty product-one subsequence, then there is a term $g$ in $T$ such that ${\mathsf{v}}_{g}(T) \geq 2 \ell(T) - n + 1$. In particular, if $n=p$ and $\ell(T) = p-1$, then all terms of $T$ are the same.  

\bigskip

\noindent {\bf Proof of Theorem \ref{EGZ-type2-length-p2}.} Let $\pi(U^{\ast}) = w \in Z(G)$ and $w \neq 1$. Without loss of generality, assume that $g_1 \bdot \dotsc \bdot g_{p-1} \bdot g_p$ is a principal part of an EGZ subsequence of $U$. For each $1 \leq i \leq p-1$, define $w_i \in Z(G)$ such that 
\[
g_p \big( g_i \dotsc g_{p-1} \big) = \big( g_i \dotsc g_{p-1} \big) g_p w_i. 
\] 
By Definition \ref{EGZ-seq}, we have $w_1, \dotsc, w_{p-1}$ are pairwise distinct and non-trivial elements of $Z(G)$, and consequently, $w_i = w^{p-1}$ for some $i$. Then we have
\[
\big( g_1 \dotsc g_{i-1} \big) g_p \big( g_i \dotsc g_{p-1} \big) g_{p+1} \dotsc g_m = 1.
\] 

\noindent For the second part, let $S$ be an EGZ sequence with $\ell(S) \geq p+1$. Then $S$ contains a subsequence $T = g_1 \bdot \dotsc \bdot g_p \bdot g_{p+1}$, where $g_1 \bdot \dotsc \bdot g_p$ is a principal part of $S$. From the previous part, it follows that $1 \in \Pi(T)$. \QED

\bigskip

\noindent The following Lemma is a consequence of the definition of the small Davenport constant and will be used in the rest of the article on several occasions:

\bigskip

\subsection{Lemma}\label{Olson-minus1} Let $G$ be a finite group and $S$ be a sequence over $G$ of length $\ell(S) = {\mathsf{d}}(G)$. If $1 \not\in \Pi(S)$, then $\Pi(S) = G \setminus \{ 1 \}$.

\bigskip

\noindent {\bf Proof.} Let $m := {\mathsf{d}}(G)$, and we write $S = g_1 \bdot \dotsc \bdot g_m$. Consider any $1 \neq g \in G$. Then, from the definition of ${\mathsf{d}}(G)$, it follows that $1 \in \Pi(S \bdot (g^{-1}))$. Let $T \mid S \bdot (g^{-1})$ such that $\pi(T^{\ast}) = 1$. Since $\ell(S) = {\mathsf{d}}(G)$, it follows that $T$ cannot be a subsequence of $S$. This implies that $\pi(T^{\ast}) = g_{i_1} \dotsc g_{i_{k-1}} g^{-1} g_{i_{k+1}} \dotsc g_{i_r}$ for distinct indices $1 \leq i_1, \dotsc, i_{k-1}, i_{k+1}, \dotsc, i_r \leq m$ with $r \geq 1$ and $1 \leq k \leq r+1$. Then $g_{i_{k+1}} \dotsc g_{i_r} g_{i_1} \dotsc g_{i_{k-1}} = g$. \QED    
    
\bigskip

\noindent For the remaining part of this section, $S$ is going to denote a sequence $g_1 \bdot \dotsc \bdot g_N$ over $G$ of length $N = 3p-2$. Set $I := \{ 1, \dotsc, N \}$, and for each $0 \leq i \leq p+1$, let $I^{S}_i$ denote the set 
\[
I^{S}_i := \{ t \in I ~:~ g_t \in {\mathcal{K}}_i \},
\]
where ${\mathcal{K}}_i$ denotes the $z$-classes as in Definition \ref{notation-z-classes}.

\bigskip

\subsection{Theorem}\label{at-least-p-2-central} (i) If $|I^{S}_0| \geq p-1$, then the sequence $S$ has a product-one subsequence. 

\smallskip

\noindent (ii) Let $|I^{S}_0| = p-2$ and denote by $U$, the subsequence of $S$ that involve the terms of $S$ from $G \setminus Z(G)$. Then $S$ contains a product-one subsequence if it satisfies one of the following three conditions:

\smallskip

\noindent (a) $|I^{S}_j| \geq p+1$ for some $j \in [1,p+1]$, \\
(b) $U$ contains a central product subsequence with at least two terms from distinct non-central $z$-classes, \\
(c) $U$ contains at least two disjoint central product subsequences.

\bigskip

\noindent {\bf Proof.} (i) If $|I^{S}_0| \geq p$, then using Olson's Theorem \ref{Ols}, we have $1 \in \Pi(S)$. Now let $|I^{S}_0| = p-1$, and without loss of generality, let $g_{2p}, \dotsc, g_N \in Z(G)$. Set $T := g_1 \bdot \dotsc \bdot g_{2p-1}$. Since ${\overline{G}} := G/{Z(G)} \cong C_p \times C_p$, from Theorem \ref{Ols}, there exists $T_1 \mid T$ with $\pi(T^{\ast}_1) \in Z(G)$. Using Olson's Theorem \ref{Ols} again, this implies that the sequence $\pi(T^{\ast}_1) \bdot g_{2p} \bdot \dotsc \bdot g_N$ (and hence $S$) has a product-one subsequence.

\smallskip

\noindent (ii) Now, assume $|I^{S}_0| = p-2$ and without loss of generality we write $U = g_1 \bdot \dotsc \bdot g_{2p}$, $g_{2p+1}, \dotsc, g_N \in Z(G)$. If at least $p+1$ terms of $U$ belong to a single $z$-class, say ${\mathcal{K}}$, then together with $g_{2p+1}, \dotsc, g_N$, we obtain a sequence of length $\geq 2p-1$ which belongs to the elementary abelian subgroup ${\mathcal{K}} \cup Z(G)$ of order $p^2$ and hence would contain a product-one subsequence. This proves (a).

\smallskip

\noindent Now let $U_1 = g_{i_1} \bdot \dotsc \cdot g_{i_{r-1}} \bdot g_{i_r}$ be a central product subsequence of $U$ with at least two terms from distinct non-central $z$-classes. Since $\pi(U_1) \subseteq Z(G)$, we may assume that $g_{i_{r-1}}, g_{i_r}$ belongs to distinct non-central $z$-classes of $G$. If the sequence $\pi(U^{\ast}_1) \bdot g_{2p+1} \bdot \dotsc \bdot g_N$ of length $p-1$ with terms from $Z(G)$ contains a product-one subsequence, then we are done. Otherwise from Lemma \ref{BEN-Lemma}, we have $\pi(U^{\ast}_1) = g_{2p+1} = \dotsc = g_N \neq 1$. Now, we have 
\[
g_{i_1} \dotsc g_{i_{r-2}} g_{i_r} g_{i_{r-1}} \neq g_{i_1} \dotsc g_{i_{r-2}} g_{i_{r-1}} g_{i_r}.
\] 
Then $g_{i_1} \dotsc g_{i_{r-2}} g_{i_r} g_{i_{r-1}} = \pi(U^{\ast}_1)^{p-\lambda}$ for some $1 \leq \lambda \leq p-2$ and consequently, \\ $g_{i_1} \dotsc g_{i_{r-2}} g_{i_r} g_{i_{r-1}} g_{2p+1} \dotsc g_{2p+\lambda} = 1$. This proves (b).

\smallskip

\noindent Next suppose $U_1$ and $U_2$ be two disjoint central product subsequence of $U$. Then, \\ $\pi(U^{\ast}_1), \pi(U^{\ast}_2), g_{2p+1}, \dotsc, g_N$ are $p$ elements of $Z(G)$ which yield a product-one subsequence of $S$. This proves (c). \QED

\bigskip

\subsection{Theorem}\label{no-central-fat-class-size} If $|I^{S}_0| = 0$ and $|I^{S}_j| \geq 2p-2$ for some $1 \leq j \leq p+1$, then $S$ has a product-one subsequence. 

\bigskip

\noindent {\bf Proof.} We fix the index $j$ with $|I^{S}_j| \geq 2p-2$ as in the statement. If $|I^{S}_j| \geq 2p-1$, then there are $2p-1$ elements in the elementary abelian normal subgroup $K := {\mathcal{K}}_j \cup Z(G)$, and we are done. So, assume that $|I^{S}_j| = 2p-2$. Without loss of generality, assume $g_1, \dotsc, g_{2p-2} \in {\mathcal{K}}_j$. Since $G/K \cong C_p$, the subsequence $g_{2p-1} \bdot \dotsc \bdot g_N$ of length $p$ must have a subsequence $T$ with product $\pi(T^{\ast}) \in K$. Then $g_1 \bdot \dotsc \bdot g_{2p-2} \bdot \pi(T^{\ast})$ will have a product-one subsequence. \QED

\bigskip

\subsection{Theorem}\label{no-central-optimum-class-size} Let $|I^{S}_0| = 0$ and $|I^{S}_j| = p$ for some $1 \leq j \leq p+1$. Assume that the terms of the sequence that belong to ${\mathcal{K}}_j$ are conjugate to each other. Then $S$ has a product-one subsequence. 

\bigskip

\noindent {\bf Proof.} Without loss of generality, assume that $g_1, \dotsc, g_p \in {\mathcal{K}}_j$. By our hypothesis, the subsequence $U := g_1 \bdot \dotsc \bdot g_p$ has central product, and no subsequence of $U$ of strictly smaller length has central product.

\smallskip

\noindent First, assume that the subsequence $T := g_{p+1} \bdot \dotsc \bdot g_N$ of length $2p-2$ has a subsequence $g_{j_1} \bdot \dotsc \bdot g_{j_r}$ with central product. Then the sequence $V$ is given by
\[
V = g_2 \bdot \dotsc \bdot g_p \bdot g_{j_1} \bdot \dotsc \bdot g_{j_r} \bdot g_1
\]
has central product. This implies that $V$ is an EGZ sequence with a principal part $g_2 \bdot \dotsc \bdot g_p \bdot g_{j_1}$. From Theorem \ref{EGZ-type2-length-p2}, it follows that $S$ has a product-one subsequence. 

\smallskip

\noindent Now, assume that $T$ has no subsequence with central product. From Lemma \ref{Olson-minus1}, there exists a subsequence $g_{j_1} \bdot \dotsc \bdot g_{j_r}$ of $T$ such that $g_p \equiv g_{j_1} \dotsc g_{j_r}$ modulo $Z(G)$. This implies that the sequence $U^{\prime}$ defined by
\[
U^{\prime} = g_1 \bdot \dotsc \bdot g_{p-1} \bdot g_{j_1} \bdot \dotsc \bdot g_{j_r}
\] 
is an EGZ sequence with principal part $g_1 \bdot \dotsc \bdot g_{p-1} \bdot g_{j_1}$. Using the Theorem \ref{EGZ-type2-length-p2} again, it follows that $S$ has a product-one subsequence. \QED

\bigskip

\noindent Before we conclude this section, we prove the following Lemma, which will be useful in the next section. 

\bigskip

\subsection{Lemma}\label{length-p3-subseqs} Let $p$ be an odd prime and $G$ be the non-abelian group of order $p^3$ and exponent $p$. Let $W$ be a sequence of length $p^3 + 3p - 3$ over $G$. Then:

\smallskip

\noindent (i) there exists a subsequence $W_1$ of $W$ of length $p^3$ whose product is central. 

\smallskip

\noindent Moreover, either $W$ has a product-one subsequence of length $p^3$ or one of the following holds:

\smallskip

\noindent (ii) at most $p^3 + p - 2$, many terms of $W$ are from the center $Z(G)$.

\smallskip

\noindent (iii) at most $p^3 + 2p - 3$, many terms of $W$ are from a maximal subgroup of $G$.

\bigskip

\noindent {\bf Proof.} (i) We have $G/{Z(G)} \cong C_p \times C_p$ and hence ${\mathsf{E}}(G/{Z(G)}) = p^2 + 2p - 2$ (Corollary 1, \cite{GAO1996}). Then there are subsequences $S_1, \dotsc, S_r$; each of them is of length $p^2$ and has central product, where $r$ can be chosen until $\sum_{i=1}^{r} \ell(S_i) < p^2 + 2p - 2$. Then we have
\[
(p^2 + 2p - 2) - p^2 \leq p^3 + 3p - 3 - rp^2 \leq p^2 + 2p - 3.
\]
This implies that
\[
(p-1) + {\frac {1}{p}} \leq r \leq p + {\frac {p-1}{p}}
\]
and consequently, $r=p$. Then $S_1 \bdot \dotsc \bdot S_p$ is a central product sequence of length $p^3$.

\smallskip

\noindent (ii) First, assume that at least $p^3 + p - 1$ terms of $S$ are central. Since ${\mathsf{E}}(Z(G)) = 2p-1$ (Corollary 1, \cite{GAO1996}), there are subsequences $S_1, \dotsc, S_r$ each of them is of length $p$, and product-one and $r$ can be chosen until $\sum_{i=1}^{r} \ell(S_i) < 2p-1$. Then we have
\[
(2p-1) - p \leq p^3 + p - 1 - rp < 2p-1.
\]
This implies that
\[
p^2 - {\frac {p-1}{p}} \leq r \leq p^2
\]
and consequently, $r = p^2$. Then $S_1 \bdot \dotsc \bdot S_{p^2}$ is a product-one sequence of length $p^3$.

\smallskip

\noindent (iii) Suppose that $p^3 + 2p - 2$ many terms of $S$ are contained in a maximal subgroup $M \leq G$. Since $M \cong C_p \times C_p$ we have ${\mathsf{E}}(M) = p^2 + 2p - 2$ (Corollary 1, \cite{GAO1996}). Then there are subsequences $S_1, \dotsc, S_r$ each of which is of length $p^2$ and product-one and $r$ can be chosen until $\sum_{i=1}^{r} \ell(S_i) < p^2 + 2p - 2$. Then we have
\[
(p^2 + 2p - 2) - p^2 \leq p^3 + 2p - 2 - rp^2 \leq p^2 + 2p - 3.
\]
This implies that
\[
(p-1) + {\frac {1}{p^2}} \leq r \leq p 
\]
and consequently, $r=p$. Then $S_1 \bdot \dotsc \bdot S_p$ is a product-one sequence of length $p^3$. \QED   

\bigskip

\section{Application to group of order $27$ and exponent $3$}
\label{Heise}

\bigskip

\noindent In this section, we will carry out the explicit details for the group $G$ of order $27$ and exponent $3$. We will use the notations from section \ref{Heisenberg_arbitrary_p}. We first prove that the small Davenport constant ${\mathsf{d}}(G)$ is equal to $6$. 

\bigskip

\subsection{Theorem}\label{Davenport-H27} Let $S := g_1 \bdot \dotsc \bdot g_7$ be a sequence of length seven over $G$. Then $S$ contains a product-one subsequence. This implies that ${\mathsf{d}}(G) = 6$.

\bigskip

\noindent {\bf Proof.} If $|I^{S}_0| \geq 2$, then the statement follows from Theorem \ref{at-least-p-2-central} (i). In the remaining part of the proof, we will adopt the following notation: we write ${\mathcal{P}}^{S}:= (|I^{S}_{j_1}|,  |I^{S}_{j_2}|, |I^{S}_{j_3}|, |I^{S}_{j_4}|) = (\mu_1, \mu_2, \mu_3, \mu_4)$ as a partition with the convention $\mu_1 \geq \mu_2 \geq \mu_3 \geq \mu_4$, where $\{ j_1, j_2, j_3, j_4 \} = [1,4]$. Without loss of generality we will assume that $g_{\mu_t + 1}, \dotsc, g_{\mu_{t+1}} \in {\mathcal{K}}_{j_{t+1}}$ for $0 \leq t \leq 3$ (with the convention $\mu_0 := 0$). We will denote by $U$ the subsequence of $S$ that consists of all non-central elements of $S$. 

\smallskip 

\noindent Now suppose that $|I^{S}_0| = 1$. Then ${\mathcal{P}}^{S}$ is a partition of $\ell(U) = 6$. If $\mu_1 \geq 4$, then the statement follows from Theorem \ref{at-least-p-2-central} (ii)(a). So we assume that $2 \leq \mu_1 \leq 3$.

\smallskip

\noindent {\bf Case I.} $\ell(U) = 6$ and $\mu_1 = 3$. 

\smallskip

\noindent Here, the only possible partitions are $(3,3,0,0), (3,2,1,0)$ and $(3,1,1,1)$.

\smallskip

\noindent If ${\mathcal{P}}^{S} = (3,3,0,0)$, then using Theorem \ref{Ols}, the subsequences $g_1 \bdot g_2 \bdot g_3$ and $g_4 \bdot g_5 \bdot g_6$ would contain two central product subsequences which are disjoint. Then, the statement follows from Theorem \ref{at-least-p-2-central} (ii)(c).

\smallskip

\noindent Now suppose ${\mathcal{P}}^{S} = (3,2,1,0)$. If $g_4 g_5 \in Z(G)$, then taking the central product subsequence of $g_1 \bdot g_2 \bdot g_3$ (which exists by Theorem \ref{Ols}), we obtain two disjoint central product subsequence of $S$ with no central terms. Using Theorem \ref{at-least-p-2-central} (ii)(c), we are done. So we assume that $g_4 g_5 \not\in Z(G)$. 

\smallskip

\noindent We claim that one of the three products $g_1 g_2, g_1 g_3, g_2 g_3$ must be non-central. If this is not true, then 
\[
(g_1 g_2)(g_1 g_3)(g_2 g_3) \equiv (g_1 g_2 g_3)^2 \hspace*{.2in} {\mathrm{mod}} \hspace*{.2in} Z(G),
\]
and consequently, $g_1 g_2 g_3 \in Z(G)$, since $G/{Z(G)}$ is elementary abelian of exponent $3$. Then $g_1 = (g_1 g_2 g_3)(g_2 g_3)^{-1} \in Z(G)$, a contradiction. This proves our claim. Now, from the claim, it follows that one of the three length $5$ sequences
\[
g_1 \bdot g_2 \bdot g_4 \bdot g_5 \bdot g_6, \hspace*{.2in} g_1 \bdot g_3 \bdot g_4 \bdot g_5 \bdot g_6, \hspace*{.2in} g_2 \bdot g_3 \bdot g_4 \bdot g_5 \bdot g_6
\]
contains a central product subsequence with at least two terms from distinct $z$-classes, and we are done by Theorem \ref{at-least-p-2-central} (ii)(b).

\smallskip

\noindent Next assume that ${\mathcal{P}}^{S} = (3,1,1,1)$. In this case, again, one of the three products $g_1 g_2, g_1 g_3, g_2 g_3$ must be non-central, and the proof is the same as above. 

\smallskip

\noindent {\bf Case II.} $\ell(U) = 6$ and $\mu_1 = 2$. 

\smallskip

\noindent Here, the only possible partitions are $(2,2,2,0), (2,2,1,1)$.

\smallskip

\noindent First assume ${\mathcal{P}}^{S} = (2,2,2,0)$. If at least two of the products $g_1 g_2, g_3 g_4, g_5 g_6$ are central, then the statement follows from Theorem \ref{at-least-p-2-central} (ii)(c). So without loss of generality, assume that either $g_1 g_2 \in Z(G)$, or else none of these products are central. Then the sequence $g_1 \bdot g_3 \bdot g_4 \bdot g_5 \bdot g_6$ of length $5$ contains a central product subsequence with at least two terms from distinct $z$-classes, and we are done by Theorem \ref{at-least-p-2-central} (ii)(b). The proof for the case ${\mathcal{P}}^{S} = (2,2,1,1)$ is again similar.

\smallskip

\noindent Now we assume that $S$ contains no central elements; i.e., $|I^{S}_{0}| = 0$. Then ${\mathcal{P}}^{S}$ is a partition of $7$. The proof will be divided into several possibilities for the partition.

\smallskip 

\noindent If $\mu_1 \geq 4$, then the proof follows from the Theorem \ref{no-central-fat-class-size}. So we assume that $2 \leq \mu_1 \leq 3$.

\smallskip

\noindent {\bf Case III.} $\ell(U) = 7$ and $\mu_1 = 3$. 

\smallskip

\noindent Here, the possible partitions are $(3,3,1,0), (3,2,2,0)$, and $(3,2,1,1)$.

\smallskip

\noindent {\bf Subcase IIIA.} The partition $(3,3,1,0)$. 

\smallskip

\noindent If the elements of $S$ that are from ${\mathcal{K}}_{j_1}$ or from ${\mathcal{K}}_{j_2}$ are pairwise conjugate to each other, then the statement follows from Theorem \ref{no-central-optimum-class-size}. Hence, without loss of generality, assume that 
\[
g_2 = g_1 u_2, g_3 = g^2_1 u_3, \hspace*{.2in} g_5 = g_4 w_5, g_6 = g^2_4 w_6
\] 
for some $u_2, u_3, w_5, w_6 \in Z(G) = \langle v \rangle$. If $u_3 = 1$ (resp. $u_2 u_3 = 1$), then $g_1 g_3 = 1$ (resp. $g_2 g_3 = 1$). Hence, we assume that the possible values of $(u_2, u_3)$ are $(1,v), (1, v^2), (v,v), (v^2, v^2)$. Using identical arguments, we assume that $(w_5, w_6)$ also assumes these values. Now the following table lists the possible product-one sequences if $(u_2, u_3) = (v,v), (v^2,v^2)$ and $(w_5,w_6) = (1,v), (1, v^2), (v,v), (v^2, v^2)$ (In case $(w_5, w_6) = (v,v), (v^2,v^2)$, the arguments are similar). The remaining cases are $(u_2, u_3), (w_5,w_6) \in \{ (1,v), (1, v^2) \}$, which are listed in the table below.  

\smallskip

\begin{center}
\begin{tabular}{|l|l|l|l|l|l|l|}
\hline
$(u_2, u_3)$ & $(w_5, w_6)$ & Extra condition & Product-one sequence \\ \hline
$(v,v)$	& $(1,v)$ & - & $g_2 \bdot g_3 \bdot g_4 \bdot g_6$ \\ \cline{2-4}
		& $(1,v^2)$ & - & $g_1 \bdot g_3 \bdot g_4 \bdot g_6$ \\  \cline{2-4}
		& $(v,v)$ & - & $g_1 \bdot g_3 \bdot g_5 \bdot g_6$ \\ \cline{2-4}
		& $(v^2,v^2)$ & - & $g_2 \bdot g_3 \bdot g_5 \bdot g_6$ \\ \hline
$(v^2,v^2)$ & $(1,v)$ & - & $g_1 \bdot g_3 \bdot g_4 \bdot g_6$ \\ \cline{2-4}	
		& $(1,v^2)$ & - & $g_2 \bdot g_3 \bdot g_4 \bdot g_6$ \\  \cline{2-4}
		& $(v,v)$ & - & $g_1 \bdot g_3 \bdot g_4 \bdot g_6$ \\  \cline{2-4}	
		& $(v^2,v^2)$ & - & $g_1 \bdot g_3 \bdot g_5 \bdot g_6$ \\  \hline
$(1,v)$ & $(1,v^2)$ & - & $g_1 \bdot g_3 \bdot g_4 \bdot g_6$ \\ \cline{2-4}
		& $(1,v)$ & $[g_4, g_1] = v^2$ & $g_1 \bdot g_4 \bdot g_3 \bdot g_6$ \\ \cline{3-4}
		&		& $[g_4, g_1] = v$ & $g_4 \bdot g_1 \bdot g_6 \bdot g_3$ \\ \hline
$(1,v^2)$ & $(1,v)$ & - & $g_1 \bdot g_3 \bdot g_4 \bdot g_6$ \\ \cline{2-4}
		& $(1,v^2)$ & $[g_4, g_1] = v^2$ & $g_4 \bdot g_1 \bdot g_6 \bdot g_3$ \\ \cline{3-4}
		&		& $[g_4, g_1] = v$ &	 $g_1 \bdot g_4 \bdot g_3 \bdot g_6$ \\ \hline	
\end{tabular}
\end{center}  

\smallskip

\noindent {\bf Subcase IIIB.} The partition $(3,2,2,0)$.

\smallskip

\noindent Here, we assume $g_1, g_2, g_3 \in {\mathcal{K}}_{j_1}, g_4, g_5 \in {\mathcal{K}}_{j_2}$ and $g_6, g_7 \in {\mathcal{K}}_{j_3}$. From Theorem \ref{Ols}, $g_1 \bdot g_2 \bdot g_3$ has a subsequence with central product. Hence, if $g_4 g_5, g_6 g_7 \in Z(G)$, then using Theorem \ref{Ols}, we are done. Next, if the terms of $g_1 \bdot g_2 \bdot g_3$ are pairwise conjugate to each other, then Theorem \ref{no-central-optimum-class-size} is applicable to produce a product-one subsequence of $S$. So, we may assume that $g_2 = g_1 u_2, g_3 = g^2_1 u_3$ for some $u_2, u_3 \in Z(G)$ with $u_3 \neq 1$ and $u_2 u_3 \neq 1$.

\smallskip

\noindent Now, we assume $g_4 g_5 \not\in Z(G)$ and $g_6 g_7 \in Z(G)$ (the case $g_4 g_5 \in Z(G)$ and $g_6 g_7 \not\in Z(G)$ is similar). Then we may assume that $g_5 = g_4 w_5$ and $g_7 = g^2_6 v_7$ for some $w_5, v_7 \in Z(G)$ with $v_7 \neq 1$. The product-one subsequences of $S$ are listed in the following table:

\smallskip

\begin{center}
\begin{tabular}{|l|l|l|l|l|l|l|}
\hline
Condition I & Condition II & Condition III & Product-one sequence \\ \hline
$u_3 \neq v_7$ & - & - & $g_1 \bdot g_3 \bdot g_6 \bdot g_7$ \\ \hline
$u_3 = v_7$ & $u_2 = v_7$ & - & $g_2 \bdot g_3 \bdot g_6 \bdot g_7$ \\ \hline
$u_3 = v_7$ & $u_2 = v^2_7$ & - & $g_2 \bdot g_3$ \\ \cline{2-4}
			& $u_2 = 1$ & $[g_6, g_1] = u_3$ & $g_3 \bdot g_6 \bdot g_1 \bdot g_7$ \\ \cline{3-4}
			&  & $[g_6, g_1] = u^2_3$ & $g_1 \bdot g_6 \bdot g_3 \bdot g_7$ \\		
\hline	
\end{tabular}
\end{center} 

\smallskip

\noindent Now, we assume $g_4 g_5, g_6 g_7 \not\in Z(G)$. Then we write $g_5 = g_4 w_5$ and $g_7 = g_6 v_7$ for some $w_5, v_7 \in Z(G)$.

\smallskip

\noindent Now, let $g_4 g_6 \in {\mathcal{K}}_{j_1}$. Then $g_4 g_6$ is conjugate to $g_5 g_7$, and consequently the terms of the sequence $g_1 \bdot g_2 \bdot g_3 \bdot g_4 g_6 \bdot g_5 g_7$ of length five belong to the elementary abelian subgroup ${\mathcal{K}}_{j_1} \cup Z(G)$. Hence, by Theorem \ref{Ols}, we have a product-one subsequence. 

\smallskip

\noindent Now, suppose $g_4 g_6 \not\in {\mathcal{K}}_{j_1}$. Then $g_4 g_6 \in {\mathcal{K}}_{j_4}$ and consequently $g^2_4 g_6 \in {\mathcal{K}}_{j_1}$, since $g_4 \not\in {\mathcal{K}}_{j_4} \cup Z(G) \leq G$. Now if $g^2_4 g_6$ is conjugate to $g_1$, then the sequence $g_1 \bdot g_2 \bdot g_4 \bdot g_5 \bdot g_6$ is an EGZ sequence with a principal part $g_1 \bdot g_2 \bdot g_4$. Next, if $g^2_4 g_6$ is conjugate to $g^2_1$, then the sequence $g_1 \bdot g_2 \bdot g_4 \bdot g_5 \bdot g_6 \bdot g_3$ is an EGZ sequence with a principal part $g_1 \bdot g_2 \bdot g_4$.

\bigskip

\noindent {\bf Subcase IIIC.} The partition $(3,2,1,1)$.

\smallskip

\noindent Here, we assume $g_1, g_2, g_3 \in {\mathcal{K}}_{j_1}, g_4, g_5 \in {\mathcal{K}}_{j_2}, g_6 \in {\mathcal{K}}_{j_3}$ and $g_7 \in {\mathcal{K}}_{j_4}$. If the terms of $g_1 \bdot g_2 \bdot g_3$ belong to the same conjugacy class, then we conclude by Theorem \ref{no-central-optimum-class-size}. Without loss of generality, we may write $g_2 = g_1 u_2, g_3 = g^2_1 u_3$ for some $u_2, u_3 \in Z(G)$. As before, we may assume $u_3 \neq 1$ and $u_2 u_3 \neq 1$.  

\smallskip

\noindent Since $g_6 g_7 \not\in {\mathcal{K}}_{j_3} \cup {\mathcal{K}}_{j_4}$, we have either $g_6 g_7 \in {\mathcal{K}}_{j_1}$ or $g_6 g_7 \in {\mathcal{K}}_{j_2}$. 

\smallskip

\noindent First suppose $g_6 g_7 \in {\mathcal{K}}_{j_2}$. Then the sequence $g_4 \bdot g_5 \bdot g_6 g_7$ has a central product subsequence. If $g_4 g_6 g_7 \in Z(G)$ (resp. $g_5 g_6 g_7 \in Z(G)$), then $g_6 \bdot g_7 \bdot g_1 \bdot g_3 \bdot g_4$ (resp. $g_6 \bdot g_7 \bdot g_1 \bdot g_3 \bdot g_5$) is an EGZ sequence with principal part $g_6 \bdot g_7 \bdot g_1$. Finally, if $g_4 g_6 g_7, g_5 g_6 g_7 \not\in Z(G)$, then $g_4, g_5$ and $g_6 g_7$ are pairwise conjugate to each other, and hence $g_4 g_5 g_6 g_7 \in Z(G)$. Then $g_4 \bdot g_5 \bdot g_1 \bdot g_3 \bdot g_6 \bdot g_7$ is an EGZ sequence with a principal part $g_4 \bdot g_5 \bdot g_1$. 

\smallskip

\noindent Now, suppose $g_6 g_7 \in {\mathcal{K}}_{j_1}$. If $g_6 g_7 = g_1 u_4$ for some $u_4 \in Z(G)$, then $g_1 \bdot g_2 \bdot g_6 \bdot g_7$ is an EGZ sequence with a principal part $g_1 \bdot g_2 \bdot g_6$. Next, if $g_6 g_7 = g^2_1 u^{\prime}_4$ for some $u^{\prime}_4 \in Z(G)$, then $g_1 \bdot g_2 \bdot g_6 \bdot g_7 \bdot g_3$ is an EGZ sequence with a principal part $g_1 \bdot g_2 \bdot g_6$. 

\smallskip

\noindent {\bf Case IV.} $\ell(U) = 6$ and $\mu_1 = 2$. 

\smallskip

\noindent Here, the only possible partition is $(2,2,2,1)$. We assume $g_1, g_2 \in {\mathcal{K}}_{j_1}, g_3, g_4 \in {\mathcal{K}}_{j_2}, g_5, g_6 \in {\mathcal{K}}_{j_3}$ and $g_7 \in {\mathcal{K}}_{j_4}$. If $g_1 g_2, g_3 g_4, g_5 g_6 \in Z(G)$, then $S$ has a product-one subsequence. So, we assume that at most two from $g_1 g_2, g_3 g_4, g_5 g_6$ are central. 

\smallskip

\noindent Now, assume that $g_1 g_2 \not\in Z(G)$ and $g_3 g_4, g_5 g_6 \in Z(G)$. 

\bigskip

\noindent Then, we have $g_3 g_5 g_4 g_6 \in Z(G)$. Suppose $g_3 g_5 g_4 g_6 = v \neq 1$. Note that $[g_5, g_3] \neq 1$. If $[g_5, g_3] = v^2$, then $g_5 g_3 g_4 g_6 = 1$. If $[g_5, g_3] = v$, then $[g_6, g_4] = v^4$, and consequently, $g_5 g_3 g_6 g_4 = 1$.

\bigskip

\noindent Next, suppose $g_1 g_2, g_3 g_4 \not\in Z(G)$. We write $g_2 = g_1 u_2, g_4 = g_3 w_4$ for some $u_2, w_4 \in Z(G)$.

\smallskip

\noindent We have $g_1 g_3 \not\in {\mathcal{K}}_{j_1} \cup {\mathcal{K}}_{j_2}$. First, suppose $g_1 g_3 \in {\mathcal{K}}_{j_3}$. Then, it follows that $g_1 \bdot g_2 \bdot g_3 \bdot g_4$ does not have a central product subsequence. From Lemma \ref{Olson-minus1}, either $g_1 g_2 g_3 g_4 g_7 \in Z(G)$ or there exist three distinct indices $1 \leq \alpha_1, \alpha_2, \alpha_3 \leq 4$ such that $g_{\alpha_1} g_{\alpha_2} g_{\alpha_3} g_7 \in Z(G)$. In the first case, $g_1 \bdot g_2 \bdot g_7 \bdot g_3 \bdot g_4$ is an EGZ sequence with a principal part $g_1 \bdot g_2 \bdot g_7$, and in the second case, $g_{\alpha_1} \bdot g_{\alpha_2} \bdot g_7 \bdot g_{\alpha_3}$ is an EGZ sequence with a principal part $g_{\alpha_1} \bdot g_{\alpha_2} \bdot g_7$.

\smallskip

\noindent Next, if $g_1 g_3 \not\in {\mathcal{K}}_{j_3}$, we have $g_1 g_3 \in {\mathcal{K}}_{j_4}$, and then $g_1 \bdot g_2 \bdot g_3 \bdot g_4 \bdot g_5$ contains an EGZ subsequence. \QED

\bigskip

\noindent The next result confirms Conjecture 2 for the non-abelian group of order $27$ and exponent $3$. The following definition is going to be helpful in the proof:

\bigskip

\subsection{Definition}\label{short-3-seq} Let $G$ be the non-abelian group of order $27$ and exponent $3$. A sequence $x_1 \bdot x_2 \bdot x_3$ over $G$ is called a {\it short $3$-sequence} if $x_1 x_2 x_3 \in Z(G)$. A short $3$-sequence is called 

\smallskip

\noindent (i) {\it central} if $x_1, x_2, x_3 \in Z(G)$, 

\smallskip

\noindent (ii) {\it thin type} if exactly two of its terms are non-central, 

\smallskip

\noindent (iii) {\it thick type A} if $x_1, x_2, x_3 \not\in Z(G)$ and the terms $x_1, x_2, x_3$ are pairwise conjugate to each other.  

\smallskip

\noindent (iv) {\it thick type B} if $x_1, x_2, x_3 \not\in Z(G)$ and it is not of thick type A. 

\smallskip

\noindent A short $3$-sequence over $G$ is said to be of {\it thick type} if it is either of thick type A or B.

\bigskip

\noindent The proof of the second part of Lemma \ref{length-p3-subseqs} applied to a cyclic group of order $3$ implies that a sequence of length $29$ always has a product-one subsequence of length $27$. We need a specific modification of this statement as follows:

\bigskip

\subsection{Lemma}\label{length28-lemma} (a) Let $T$ be a sequence of length $28$ over the cyclic group $C_3$ of order $3$. Then either $T$ has a product-one subsequence of length $27$, or it has subsequences $T^{\prime}$ and $T^{\prime \prime}$ of length $27$ such that $\pi((T^{\prime})^{\ast})$ and $\pi((T^{\prime \prime})^{\ast})$ are distinct and both non-trivial.

\smallskip

\noindent (b) Let $T$ be a sequence of length $27$ over the cyclic group $C_3 = \{ 1, v, v^2 \}$ of order $3$. If $T$ is not product-one, then it has subsequences $U_1, \dotsc, U_8$ each of length $3$ and containing the same terms, and the remaining three terms form a subsequence of one of the following forms (up to a permutation):
\[
1 \bdot 1 \bdot v, \hspace*{.2in} 1 \bdot 1  \bdot v^2, \hspace*{.2in} 1 \bdot v \bdot v, \hspace*{.2in} 1 \bdot v^2 \bdot v^2, \hspace*{.2in} v \bdot v \bdot v^2, \hspace*{.2in} v \bdot v^2 \bdot v^2.
\] 

\bigskip

\noindent {\bf Proof.} (a) We write $C_3 = \langle u \rangle$. Suppose $T$ has no product-one subsequence of length $27$. Then there are sequences $T_1, \dotsc, T_8$ each of length $3$, and each $T_i$ is of the form $T_i = v \bdot v \bdot v$ for some $v \in C_3$. By our hypothesis, the remaining four terms must be of the form $(\alpha, \alpha, \beta, \beta)$ for some distinct $\alpha, \beta \in C_3$. Based on these choices, we define the sequences $T^{\prime}_9$ and $T^{\prime \prime}_9$, which are listed below.

\begin{center}
\begin{tabular}{|l|l|l|l|l|l|l|}
\hline
$(\alpha, \alpha)$ & $(\beta, \beta)$ & $T^{\prime}_9$ & $T^{\prime \prime}_9$ \\ \hline
$(1,1)$ & $(u,u)$ & $1 \bdot 1 \bdot u$ & $1 \bdot u \bdot u$ \\ \hline
$(1,1)$	& $(u^2,u^2)$ & $1 \bdot u^2 \bdot u^2$ & $1 \bdot 1 \bdot u^2$ \\ \hline
$(u,u)$ & $(u^2,u^2)$ & $u \bdot u \bdot u^2$ & $u \bdot u^2 \bdot u^2$ \\ \hline
\end{tabular}
\end{center} 

\noindent Finally, we define 
\[
T^{\prime} := T_1  \bdot \dotsc \bdot T_8  \bdot T^{\prime}_9, \hspace*{.2in} T^{\prime \prime} := T_1  \bdot \dotsc \bdot T_8 \bdot T^{\prime \prime}_9.
\] 
Notice that $\pi((T^{\prime})^{\ast}) = u$ and $\pi((T^{\prime \prime})^{\ast}) = u^2$. 

\bigskip

\noindent (b) As shown in part (a), we can choose subsequences $T_i ~(1 \leq i \leq 7)$ that have the desired form. Let $T^{\prime} := T \bdot (T_1 \bdot \dotsc \bdot T_7)^{[-1]}$. If $T^{\prime}$ has no element with multiplicity $\geq 3$, then each element of $C_3$ must occur exactly twice among them, and consequently $\pi((T^{\prime})^{\ast}) = 1$ and hence $\pi(T^{\ast}) = 1$. Hence, there is a subsequence $T_8$ of $T^{\prime}$ that is also of the desired form. It is easy to see that the remaining three elements have the desired form as listed if $\pi(T^{\ast}) \neq 1$. \QED

\bigskip

\subsection{Lemma}\label{twin-thick-type} Let $G$ be the non-abelian group of order $27$ and exponent $3$. Let $W, U_1$ and $U_2$ be central product sequences over $G$, where $U_1$ and $U_2$ are short $3$-sequences of non-central type. Then the following holds:

\smallskip

\noindent (a) If $W$ contains terms from four distinct $z$-classes, then $W$ is an EGZ sequence.

\smallskip

\noindent (b) If $U_1 \bdot U_2$ contains terms from at least two distinct $z$-classes and one of $U_1$ and $U_2$ is of thick type, then $U_1 \bdot U_2$ is an EGZ sequence.

\smallskip

\noindent (c) Let $g_1, g_2, g_3, g_4$ be non-central elements and $u \in Z(G)$. Suppose $g_1$ and $g_3$ are from two distinct $z$-classes and $g_1 g_2, g_3 g_4 \in Z(G)$. Then there exists a permutation $\sigma$ of $[1,4]$ such that $g_{\sigma(1)} g_{\sigma(2)} g_{\sigma(3)} g_{\sigma(4)} u = 1$. 

\bigskip

\noindent {\bf Proof.} We write $U_i = g_{i1} \bdot g_{i2} \bdot g_{i3}$ for $i=1,2$ and $U := U_1 \bdot U_2$. 

\smallskip

\noindent (a) If $W$ contains terms from four distinct non-central $z$-classes of $G$, say $y_1, y_2, y_3, y_4$, then $y_1 y_2$ belongs to the $z$-class of either $y_3$ or $y_4$. Then $U$ contains an EGZ subsequence with a principal part $y_1 \bdot y_2 \bdot y_4$ (resp. $y_1 \bdot y_2 \bdot y_3$).

\smallskip

\noindent (b) Let $U_1$ be of thick type B and $g_{1j} \in {\mathcal{K}}_j, ~j=1,2,3$. Then the $z$-classes ${\mathcal{K}}_1, {\mathcal{K}}_2, {\mathcal{K}}_3$ are all distinct. If $U_2$ is of thin type, with a permutation if necessary, we may assume that $g_{21}$ is conjugate to $g_{11}$. Then $U$ has an EGZ subsequence with a principal part $g_{11} \bdot g_{21} \bdot g_{13}$. Next, suppose that $U_2$ is of thick type A, and we may assume that elements of $U_2$ belong to the same $z$-class as $g_{11}$. Then, either $g_{11}$ is conjugate to $g_{21}$, or else $g_{11} g_{21} \in Z(G)$. Then, $U$ has an EGZ subsequence with a principal part $g_{11} \bdot g_{21} \bdot g_{13}$ (resp. $g_{21} \bdot g_{22} \bdot g_{13}$). 

\smallskip

\noindent So, we assume that both $U_1$ and $U_2$ are of thick type B, and from part (a), we assume that $g_{1j}$ and $g_{2j}$ belong to the same $z$-class. If $g_{11}$ is conjugate to $g_{21}$, then $g_{11} \bdot g_{21} \bdot g_{13}$ is a principal part of an EGZ subsequence of $U$. If $g_{11}$ is not conjugate to $g_{21}$, then we have $g_{11} g_{21}, g_{12} g_{22}, g_{13} g_{23} \in Z(G)$, and then $g_{11} \bdot g_{22} \bdot g_{13}$ is a principal part of an EGZ subsequence of $U$. 

\smallskip

\noindent Now, assume that $U_1$ is of thick type A. Since $U_1 \bdot U_2$ contains elements from at least two distinct $z$-classes, $U_2$ has a term $v$ that belongs to a $z$-class other than the terms of $U_1$. Then $g_{11} \bdot g_{12} \bdot v$ is a principal part of an EGZ subsequence of $U$. 



\smallskip

\noindent (c) By hypothesis, $g_1 = g^2_2 w$ for some $w \in Z(G)$. We assume that $g_1 g_2 g_3 g_4 z = v \neq 1$. If $[g_2, g_3] = v$, then $g_1 g_3 g_2 g_4 z = 1$, and if $[g_2, g_3] = v^2$, then $[g_1, g_3] = v$, and hence $g_2 g_3 g_1 g_4 z = 1$. \QED

\smallskip

\subsection{Theorem}\label{EGZ-H27} Let $S := g_1 \bdot \dotsc \bdot g_{33}$ be a sequence of length $33$ over $G$. Then $S$ contains a product-one subsequence of length $27$. This implies that ${\mathsf{E}}(G) = 33$.

\bigskip

\noindent{\bf Proof.} Let $n_0$ denote the number of terms of $S$ from $Z(G)$. From Lemma \ref{length-p3-subseqs}(ii), we may assume that $n_0 \leq 28$. Since ${\mathsf{s}}(G/{Z(G)}) = 9$ (see Theorem 5.8.3, \cite{GHK2006}), there exist short $3$-subsequences $U_1, \dotsc, U_9$ of $S$. We set $U := U_1 \bdot \dotsc \bdot U_9$. Next, from Lemma \ref{length-p3-subseqs}(iii), it follows that elements of $S$ are from at least two distinct non-central $z$-classes of $G$. For the rest of the proof, we write $U_i = g_{i1} \bdot g_{i2} \bdot g_{i3}$ for $i \in [1,9]$ and we divide the proof into several cases:

\smallskip

\noindent {\bf Case I.} $U$ contain elements from non-central $z$-classes of $G$.

\smallskip

\noindent If $U$ contains elements from four distinct non-central $z$-classes of $G$, then from Lemma \ref{twin-thick-type}(a) and Theorem \ref{EGZ-type2-length-p2}, we are done. So, we assume that the elements of $U$ are from at most three distinct non-central $z$-classes: ${\mathcal{K}}_1, {\mathcal{K}}_2, {\mathcal{K}}_3$.

\smallskip

\noindent {\bf Subcase IA.} At least one of the $U_i$ is of thick type B.

\smallskip

\noindent Without loss of generality, we assume that $U_1$ is of thick type B. 

\smallskip

\noindent Now if $U_j$ contain non-central elements for some $2 \leq j \leq 9$, then from Lemma \ref{twin-thick-type}(b) we have $U_1 \bdot U_j$ is an EGZ sequence and hence from Theorem \ref{EGZ-type2-length-p2}, we are done. Hence, we may assume that the short $3$-sequences $U_2, \dotsc, U_9$ are all central. 

\smallskip

\noindent Let $S \bdot U^{[-1]}$ has a subsequence $Z := z_1 \bdot z_2 \bdot z_3$ consisting of central elements. Since $U_1$ contain terms from distinct $z$-classes, up to a permutation within $U_1$ if necessary, we assume that $\pi(U^{\ast}_1) = v \neq 1$. Now we have a sequence $W := U_2 \bdot \dotsc \bdot U_9 \bdot Z$ of length $27$ over $G$, all whose terms are central. Then either $\pi(W^{\ast}) = 1$, or from Lemma \ref{length28-lemma}(b) there exists central short $3$-sequences $T_1, \dotsc, T_8 \mid W$ each of which containing same terms within itself and the remaining three central elements form a short $3$-sequence $T_9$ with $\pi(T^{\ast}_9) = v$ or $v^2$. If $\pi(T^{\ast}_9) = v^2$, then $U_1 \bdot T_2 \bdot \dotsc \bdot T_8 \bdot T_9$ is a product-one sequence of length $27$. The same proof works if $\pi(T^{\ast}_9) = v$ and $\pi((U^{\prime}_1)^{\ast}) = v^2$ for another permutation $U^{\prime}_1$ of $U_1$. Since $\lvert \pi(U_1) \rvert \geq 2$, we may assume that $\pi(U_1) = \{ 1, v \}$. In this case, we have a permutation $U^{\prime \prime}_1$ of $U_1$ such that $\pi((U^{\prime \prime}_1)^{\ast}) = 1$. Then $U^{\prime \prime}_1 \bdot T_1 \bdot \dotsc T_8$ is a product-one sequence of length $27$. 

\smallskip

\noindent This reduces to the case that at most two terms of $S \bdot U^{[-1]}$ are central (and hence at least four terms are non-central). We write $S \bdot U^{[-1]} = h_1 \bdot \dotsc \bdot h_5 \bdot h_6$. Suppose (up to a reordering) $h_1, h_2, h_3, h_4 \not\in Z(G)$ and set $W := h_1 \bdot h_2 \bdot h_3 \bdot h_4$. Then $U_1 \bdot W$ has length $7$ and consists of all non-central terms. 

\smallskip

\noindent We follow the notations of partitions introduced in the proof of Theorem \ref{Davenport-H27}. We rewrite the $z$-class indices to denote them by ${\mathcal{K}}_{j_1}, {\mathcal{K}}_{j_2}, {\mathcal{K}}_{j_3}, {\mathcal{K}}_{j_4}$ so that ${\mathcal{P}}^{U_1 \bdot W} = (\mu_1, \mu_2, \mu_3, \mu_4)$ and $\mu_1 \geq \mu_2 \geq \mu_3 \geq \mu_4$. Also, when $\mu_4 = 0$, upto a reordering if necessary, we assume $g_{11} \in {\mathcal{K}}_{j_1}, g_{12} \in {\mathcal{K}}_{j_2}$ and $g_{13} \in {\mathcal{K}}_{j_3}$. Since three entries of ${\mathcal{P}}^{U_1 \bdot W}$ are non-zero, we have $\mu_1 \leq 5$. In the following, if we obtain three terms of $W$ that belong to the same non-central $z$-class, say these be $h_1, h_2, h_3$, then either $h_1 h_2 h_3 \in Z(G)$, or else $h_1 h_2 \in Z(G)$ (up to a permutation). Then we define $W_1 := h_1 \bdot h_2 \bdot h_3$ (resp. $W_1 := h_1 \bdot h_2 \bdot g_{23}$) and consequently $U_1 \bdot W_1 \bdot U_3 \bdot \dotsc \bdot U_9$ is an EGZ sequence from Lemma \ref{twin-thick-type}(b). Then we are done by Theorem \ref{EGZ-type2-length-p2}. If $\mu_1 = 5$, then ${\mathcal{P}}^{U_1 \bdot W} = (5,1,1,0)$. Then all the four terms of $W$ belong to the same non-central $z$-class, and we are done. Next if $\mu_1 = 4$, we have ${\mathcal{P}}^{U_1 \bdot W} = (4,2,1,0)$ or $(4,1,1,1)$. In either case, we have at least three terms of $W$ belonging to the same $z$-class. 

\smallskip

\noindent So assume $\mu_1 =3$, and then ${\mathcal{P}}^{U_1 \bdot W}$ is $(3,3,1,0), (3,2,2,0)$, or $(3,2,1,1)$. 

\smallskip

\noindent If ${\mathcal{P}}^{U_1 \bdot W} = (3,3,1,0)$, then we may assume that $h_1, h_2 \in {\mathcal{K}}_{j_1}$ and $h_3, h_4 \in {\mathcal{K}}_{j_2}$. If either $h_1 h_2 \in Z(G)$ or $h_3 h_4 \in Z(G)$, then we define $W_1 := h_1 \bdot h_2 \bdot g_{23}$ (resp $W_1 := h_3 \bdot h_4 \bdot g_{23}$) and $U_1 \bdot W_1 \bdot U_3 \bdot \dotsc U_9$ is an EGZ sequence of length $27$. So we assume that $h_1$ (resp. $h_3$) is conjugate to $h_2$ (resp. $h_4$). Now if $h_1$ is conjugate to $g_{11}$ (resp. $h_1 g_{11} \in Z(G)$) and $h_3 g_{12} \in Z(G)$ (resp. $h_3$ is conjugate to $g_{12}$), then $(g_{11} \bdot h_1 \bdot h_2) \bdot (g_{12} \bdot h_3 \bdot g_{23}) \bdot U_3 \bdot \dotsc \bdot U_9$ (resp. $(g_{11} \bdot h_1 \bdot g_{23}) \bdot (g_{12} \bdot h_3 \bdot h_4) \bdot U_3 \bdot \dotsc \bdot U_9$) is an EGZ sequence of length $27$ from Lemma \ref{twin-thick-type}(b). Finally we consider the case $h_1$ is conjugate to $g_{11}$ and $h_3$ is conjugate to $g_{12}$. Then $(g_{11} \bdot h_1 \bdot h_2) \bdot (g_{12} \bdot h_3 \bdot h_4) \bdot U_3 \bdot \dotsc \bdot U_9$ is again an EGZ sequence of length $27$ from Lemma \ref{twin-thick-type}(b). 

\smallskip

\noindent Next, let ${\mathcal{P}}^{U_1 \bdot W} = (3,2,2,0)$. Here we may assume $h_1, h_2 \in {\mathcal{K}}_{j_1}, h_3 \in {\mathcal{K}}_{j_2}$ and $h_4 \in {\mathcal{K}}_{j_3}$ (up to a reordering). Now if one of the products $h_1 h_2, h_1 h_3 h_4, h_2 h_3 h_4$ is central, then we define $W_2 := h_1 \bdot h_2 \bdot g_{23}$ (resp. $h_1 \bdot h_3 \bdot h_4$ or $h_2 \bdot h_3 \bdot h_4$) and then $U_1 \bdot W_2 \bdot U_3 \bdot \dotsc \bdot U_9$ is a central product sequence and the statement follows from Lemma \ref{twin-thick-type}(b). So we assume that $h_1$ is conjugate to $h_2$. Now if $g_{11}$ is non-conjugate to $h_1$, then $g_{11}$ is conjugate to $h_1 h_2$. Then $(h_1 \bdot h_2 \bdot g_{12} \bdot g_{13} \bdot g_{22} \bdot g_{23}) \bdot U_3 \bdot \dotsc \bdot U_9$ is a EGZ sequence of length $27$ with principal part $h_1 \bdot h_2 \bdot g_{12}$. So we assume that $g_{11}, h_1, h_2$ are pairwise conjugate to each other. Now if $g_{12}$ is non-conjugate to $h_3$, then we have a length $27$ sequence
\[
(g_{11} \bdot h_1 \bdot h_2) \bdot (g_{12} \bdot h_3 \bdot g_{23}) \bdot U_3 \bdot \dotsc \bdot U_9
\]
and we are done by Lemma \ref{twin-thick-type}(b). So we assume that $g_{12}$ is conjugate to $h_3$ and similarly, $g_{13}$ is conjugate to $h_4$. Then $h_1 h_3 h_4 \in Z(G)$ and Lemma \ref{twin-thick-type}(b) applied to the sequence $U_1 \bdot (h_1 \bdot h_3 \bdot h_4) \bdot U_3 \bdot \dotsc \bdot U_9$, we are done.

\smallskip

\noindent Now we consider the case ${\mathcal{P}}^{U_1 \bdot W} = (3,2,1,1)$. In case $(g_{11}, g_{12}, g_{13}) \in {\mathcal{K}}_{j_2} \times {\mathcal{K}}_{j_3} \times {\mathcal{K}}_{j_4}$, then three elements of $W$ belongs to ${\mathcal{K}}_{j_1}$. Then up to a reordering, we may assume either $h_1 h_2 h_3 \in Z(G)$, or $h_1 h_2 \in Z(G)$. Then applying Lemma \ref{twin-thick-type}(b) to the sequence $U_1 \bdot (h_1 \bdot h_2 \bdot h_3) \bdot U_3 \bdot \dotsc \bdot U_9$ (resp. $U_1 \bdot (h_1 \bdot h_2 \bdot g_{23}) \bdot U_3 \bdot \dotsc \bdot U_9$) of length $27$ we are done. Then without loss generality, we assume that $g_{11}, h_1, h_2 \in {\mathcal{K}}_{j_1}$ and there are three possible cases: $(g_{12}, g_{13}) \in {\mathcal{K}}_{j_2} \times {\mathcal{K}}_{j_3}$ or ${\mathcal{K}}_{j_2} \times {\mathcal{K}}_{j_4}$ or ${\mathcal{K}}_{j_3} \times {\mathcal{K}}_{j_4}$. Also, from above, we assume that $h_1$ is conjugate to $h_2$. Now if $g_{11}$ is non-conjugate to $h_1$, then $(h_1 \bdot h_2 \bdot g_{12} \bdot g_{13} \bdot g_{22} \bdot g_{23}) \bdot U_3 \bdot \dotsc \bdot U_9$ is an EGZ sequence of length $27$ with principal part $h_1 \bdot h_2 \bdot g_{12}$. So we assume that $g_{11}, h_1, h_2$ are pairwise conjugate to each other. 
 
\smallskip

\noindent Now suppose $(g_{12}, g_{13}) \in {\mathcal{K}}_{j_2} \times {\mathcal{K}}_{j_3}$ and then we may assume that $h_3 \in {\mathcal{K}}_{j_2}, h_4 \in {\mathcal{K}}_{j_4}$. If $g_{12} h_3 \in Z(G)$, then Lemma \ref{twin-thick-type}(b) applied to $(g_{11} \bdot h_1 \bdot h_2) \bdot (g_{12} \bdot h_3 \bdot g_{23}) \bdot U_3 \bdot \dotsc \bdot U_9$ would work. So we assume that $g_{12}$ is conjugate to $h_3$. From Lemma \ref{Olson-minus1}, either $W$ has a central product subsequence $W_3$ (and then $3 \leq \ell(W_3) \leq 4$) or else each non-central element of $G$ is congruent to an element of $\Pi(W)$ modulo $Z(G)$. If such $W_3$ exists with $\ell(W_3) = 3$, then using Lemma \ref{twin-thick-type}(b) to $U_1 \bdot W_3 \bdot U_3 \bdot \dotsc \bdot U_9$, we are done. In case $\ell(W_3) = 4$, then 
\[
(h_1 \bdot h_2 \bdot h_3 \bdot h_4 \bdot g_{22} \bdot g_{23}) \bdot U_3 \bdot \dotsc \bdot U_9
\]
is an EGZ sequence with principal part $h_1 \bdot h_2 \bdot h_3$. Now assume that $W_3$ does not exist. Then $g^{-1}_{12} \equiv \pi(W^{\ast}_4)$ modulo $Z(G)$ for some $W_4 \mid W$ with $2 \leq \ell(W_4) \leq 4$. If $\ell(W_4) = 2$, then using multiplication properties of $z$-classes $W_4 \neq h_1 \bdot h_2$ and $h_3$ cannot appear in $W_4$ since $g_{12}$ is conjugate to $h_3$. Then $g^{-1}_{12} \equiv h_1 h_4$ and also $g^{-1}_{12} \equiv h_2 h_4$ modulo $Z(G)$. Then we have
\[
g_{11} g_{12} g_{13} \equiv 1 \equiv h_1 g_{12} h_4 \equiv g_{11} g_{12} h_4 \hspace*{.2in} {\mathrm{mod}}~Z(G),
\]
which implies that $g_{13}$ and $h_4$ belongs to the same $z$-class, a contradiction. If $\ell(W_4) = 3$ and $h_1 \bdot h_2 \mid W_4$, then $(h_1 \bdot h_2 \bdot g_{12} \bdot \xi \bdot g_{22} \bdot g_{23}) \bdot U_3 \bdot \dotsc \bdot U_9$ is an EGZ sequence of length $27$ with principal part $h_1 \bdot h_2 \bdot g_{12}$, where $\xi \in \{ h_3, h_4 \}$. Next if $\ell(W_4) = 3$ and only one of $h_1, h_2$ appears in $W_4$, then $(h_1 \bdot g_{12} \bdot h_4 \bdot h_3 \bdot g_{22} \bdot g_{23}) \bdot U_3 \bdot \dotsc \bdot U_9$ is an EGZ sequence of length $27$ with principal part $h_1 \bdot g_{12} \bdot h_4$, since if $h_1 g_{12}$ centralize $h_4$ then $h_1 g_{12} h_4 \in {\mathcal{K}}[h_4]$, which contradicts ${\mathcal{K}}[h_4] \neq {\mathcal{K}}[h_3]$. Finally, if $\ell(W_4) = 4$, then $g_{12} \bdot W_4 \bdot g_{23} \bdot U_3 \bdot \dotsc \bdot U_9$ is an EGZ sequence of length $27$ with principal part $h_1 \bdot h_2 \bdot g_{12}$. 

\smallskip

\noindent The proof for $(g_{12}, g_{13}) \in {\mathcal{K}}_{j_2} \times {\mathcal{K}}_{j_4}$ is same as in the above paragraph. So we consider the case $(g_{12}, g_{13}) \in {\mathcal{K}}_{j_3} \times {\mathcal{K}}_{j_4}$ and again we may assume that $h_1, h_2 \in {\mathcal{K}}_{j_1}, h_3, h_4 \in {\mathcal{K}}_{j_2}$ with $h_1$ is conjugate to $h_2$ and $h_3$ is conjugate to $h_4$. Next if $g_{11}$ is non-conjugate to $h_1$, then $g_{11}$ is conjugate to $h_1 h_2$ and then $(h_1 \bdot h_2 \bdot g_{12} \bdot g_{13} \bdot g_{22} \bdot g_{23}) \bdot U_3 \bdot \dotsc \bdot U_9$ is an EGZ sequence of length $27$ with principal part $h_1 \bdot h_2 \bdot g_{12}$. So we assume that $g_{11}$ is conjugate to $h_1$. Now notice that with this hypothesis, $W$ does not contain any central product subsequence, and hence by Lemma \ref{Olson-minus1}, $g^{-1}_{12} \equiv \pi(W^{\ast}_5)$ modulo $Z(G)$ for some $W_5 \mid W$. In case $\ell(W_5) = 2$, then the only possible choice is $W_5 = h_1 \bdot h_3$ (and also at the same time $h_1 \bdot h_4, h_2 \bdot h_3, h_2 \bdot h_4$ using the conjugacy hypothesis). But then
\[
g_{11} g_{12} g_{13} \equiv 1 \equiv h_1 g_{12} h_3 \equiv g_{11} g_{12} h_3 \hspace*{.2in} {\mathrm{mod}}~Z(G),
\] 
which implies that $g_{13}$ is conjugate to $h_3$, a contradiction. This implies that $3 \leq \ell(W_5) \leq 4$, and from the conjugacy hypothesis on the terms of $W$, it follows that $g_{12} \bdot W_5 \bdot g_{22} \bdot g_{23} \bdot U_3 \bdot \dotsc \bdot U_9$ (resp. $g_{12} \bdot W_5 \bdot g_{22} \bdot U_3 \bdot \dotsc \bdot U_9$) in case $\ell(W_5) = 3$ (resp. $\ell(W_5) = 4$) are EGZ sequences of length $27$. 

\smallskip

\noindent Now assume that $\mu_1 = 2$ and then ${\mathcal{P}}^{U_1 \bdot W} = (2,2,2,1)$. 

\smallskip

\noindent First assume that $(g_{11}, g_{12}, g_{13})$, $(h_1, h_2, h_3) \in {\mathcal{K}}_{j_1} \times {\mathcal{K}}_{j_2} \times {\mathcal{K}}_{j_3}$ and $h_4 \in {\mathcal{K}}_{j_4}$, up to a reordering. If $W$ has a central product subsequence $W_6 \mid W$, then $3 \leq \ell(W_6) \leq 4$ since all the terms of $W$ belong to different non-central $z$-classes. If $\ell(W_6) = 3$, then applying Lemma \ref{twin-thick-type}(b) to $U_1 \bdot W_6 \bdot U_3 \bdot \dotsc \bdot U_9$ we are done. So suppose $\ell(W_6) = 4$, which means $h_1 h_2 h_3 h_4 \in Z(G)$. Now if $h_3$ centralize $h_1 h_2$, then either $h_1 h_2 h_3 \in Z(G)$, or else $h_1 h_2$ is conjugate to $h_3$. The first case is not possible since $h_4 \not\in Z(G)$, and in the second case $h_1 h_2 h_3 \in {\mathcal{K}}[h_3] = {\mathcal{K}}_{j_3}$, which implies that $h_4 \in {\mathcal{K}}_{j_3}$ and both of these cannot be possible. This implies that $h_1 \bdot h_2 \bdot h_3 \bdot h_4 \bdot g_{22} \bdot g_{23} \bdot U_3 \bdot \dotsc \bdot U_9$ is an EGZ sequence with principal part $h_1 \bdot h_2 \bdot h_3$. Now, we assume that $W$ does not have a central product subsequence. This also implies that at least one of the tuples $(g_{11}, h_1), (g_{12}, h_2), (g_{13}, h_3)$ have non-conjugate entries (or else $h_1 h_2 h_3 \in Z(G)$) and at least one would have entries that are conjugate to each other (or else $(h_1 h_2 h_3)^2 \in Z(G)$ and hence $h_1 h_2 h_3 \in Z(G)$). Without loss of generality assume that $g_{11}$ is conjugate to $h_1$ and $g_{13} h_3 \in Z(G)$. From Lemma \ref{Olson-minus1}, $g^{-1}_{11} \equiv \pi((W^{\prime}_6)^{\ast})$ modulo $Z(G)$ for some $W^{\prime}_6 \mid W$. From above hypothesis, we necessarily would have $h_1$ is a term of $W^{\prime}_6$, or else $g_{11} \pi((W^{\prime}_6)^{\ast}) \equiv h_1 \pi((W^{\prime}_6)^{\ast})$ modulo $Z(G)$, which implies that $W$ have a central product subsequence, a contradiction. Next since $g_{11} h_1 \in {\mathcal{K}}[h_1]$, we must have $3 \leq \ell(W^{\prime}_6) \leq 4$. If $\ell(W^{\prime}_6) = 3$, we write $W^{\prime}_6 = h_1 \bdot \xi_1 \bdot \xi_2$ for some $\xi_1, \xi_2 \in \{ h_2, h_3, h_4 \}$. Then $g_{11} \bdot W^{\prime}_6 \bdot g_{22} \bdot g_{23} \bdot U_3 \bdot \dotsc \bdot U_9$ is an EGZ sequence with principal part $g_{11} \bdot h_1 \bdot \xi_1$. If $\ell(W^{\prime}_6) = 4$, then $g_{11} \bdot W^{\prime}_6 \bdot g_{23} \bdot U_3 \bdot \dotsc \bdot U_9$ is an EGZ sequence with principal part $g_{11} \bdot h_1 \bdot h_2$. 

\smallskip

\noindent The only remaining case for the above partition is when $(g_{11}, g_{12}), (h_1, h_2) \in {\mathcal{K}}_{j_1} \times {\mathcal{K}}_{j_2}$, $g_{13} \in {\mathcal{K}}_{j_4}$ and $h_3, h_4 \in {\mathcal{K}}_{j_3}$, up to a reordering. If $h_3 h_4 \in Z(G)$, then applying Lemma \ref{twin-thick-type}(b) to $U_1 \bdot (h_3 \bdot h_4 \bdot g_{23}) \bdot U_3 \bdot \dotsc \bdot U_9$, we are done. So we assume that $h_3$ is conjugate to $h_4$. Now if there is a central product subsequence $W_7 \mid W$, then $3 \leq \ell(W_7) \leq 4$. If $\ell(W_7) = 3$, then applying Lemma \ref{twin-thick-type}(b) again to $U_1 \bdot W_7 \bdot U_3 \bdot \dotsc \bdot U_9$, we are done. If $\ell(W_7) = 4$, then $h_3 \bdot h_4 \bdot h_1 \bdot h_2 \bdot g_{22} \bdot g_{23} \bdot U_3 \bdot \dotsc \bdot U_9$ is an EGZ sequence of length $27$ with principal part $h_3 \bdot h_4 \bdot h_1$. 

\smallskip

\noindent Now, we assume that $W$ does not contain a central product subsequence. If any two of the products $g_{11} h_1, g_{12} h_2, h_3 h_4$ is central, denoting them by $\xi_1 \xi_2, \xi_3 \xi_4$, we are done by applying Lemma \ref{twin-thick-type}(c) to the sequence $\xi_1 \bdot \xi_2 \bdot \xi_3 \bdot \xi_4 \bdot g_{22} \bdot g_{23} \bdot U_3 \bdot \dotsc \bdot U_9$. Next if $h_3 h_4 \in Z(G)$, then we are done by applying Lemma \ref{twin-thick-type}(b) to the sequence $U_1 \bdot (h_3 \bdot h_4 \bdot g_{23}) \bdot U_3 \bdot \dotsc \bdot U_9$. So we assume that among the products $g_{11} h_1, g_{12} h_2, h_3 h_4$, either one of the first two is central, or none of them is central. First suppose that $g_{12} h_2 \in Z(G)$ and $g_{11} h_1, h_3 h_4 \not\in Z(G)$. Then from Lemma \ref{Olson-minus1}, $g^{-1}_{11} \equiv \pi((W^{\prime}_7)^{\ast})$ modulo $Z(G)$ for some $W^{\prime}_7 \mid W$ with $1 \leq \ell(W^{\prime}_7) \leq 4$. If $\ell(W^{\prime}_7) = 1$, then $W^{\prime}_7 = h_1$, which implies that $g_{11} h_1 \in Z(G)$, a contradiction. If $\ell(W^{\prime}_7) = 2$, then $W^{\prime}_7$ cannot have the term $h_1$, since $g_{11} h_1 \in {\mathcal{K}}[h_1]$. Hence $g_{11} \xi_1 \xi_2 \in Z(G)$ for some $\xi_1, \xi_2 \in \{ h_2, h_3, h_4 \}$. But since $g_{11}$ is conjugate to $h_1$, this implies that $W$ has a central product subsequence of length $3$, a contradiction. If $\ell(W^{\prime}_7) = 3$, then again $W^{\prime}_7 \neq h_2 \bdot h_3 \bdot h_4$, using the same reason as above. Hence $W^{\prime}_7$ has one of the three possibilities $h_1 \bdot h_2 \bdot h_3, h_1 \bdot h_2 \bdot h_4, h_1 \bdot h_3 \bdot h_4$. Then we obtain an EGZ sequence of length $27$ in each of them given by $g_{11} \bdot W^{\prime}_7 \bdot g_{22} \bdot g_{23} \bdot U_3 \bdot \dotsc \bdot U_9$ with principal part $g_{11} \bdot h_1 \bdot \xi$ where $\xi \in \{ h_2, h_3 \}$. If  $\ell(W^{\prime}_7) = 4$, then $g_{11} \bdot W \bdot g_{23} \bdot U_3 \bdot \dotsc \bdot U_9$ is an EGZ sequence of length $27$ with principal part $g_{11} \bdot h_1 \bdot h_2$. Next, the case when $g_{11} h_1 \in Z(G)$ and $g_{12} h_2, h_3 h_4 \not\in Z(G)$ is similarly derived by analyzing $g^{-1}_{12} \equiv \pi((W^{\prime}_7)^{\ast})$ modulo $Z(G)$ for some $W^{\prime}_7 \mid W$ with $1 \leq \ell(W^{\prime}_7) \leq 4$. So we assume that $g_{11} h_1, g_{12} h_2, h_3 h_4 \not\in Z(G)$. The proof for this case again similar by considering either $g^{-1}_{11} \equiv \pi((W^{\prime}_7)^{\ast})$ or $g^{-1}_{12} \equiv \pi((W^{\prime}_7)^{\ast})$ modulo $Z(G)$ for some $W^{\prime}_7 \mid W$ with $1 \leq \ell(W^{\prime}_7) \leq 4$.

\bigskip

\noindent {\bf Subcase IB.} At least one of the $U_i$ is of thick type, but none are of thick type B.

\smallskip

\noindent As before, we assume that $U_1$ is of thick type A, the terms of $U \bdot U^{[-1]}_1$ are all central, the elements of $U_1$ are from ${\mathcal{K}}_1$, and the remaining distinct non-central classes are ${\mathcal{K}}_2, {\mathcal{K}}_3, {\mathcal{K}}_4$. From Lemma \ref{length-p3-subseqs}(iii), at most three elements of $S \bdot U^{[-1]}$ would be from ${\mathcal{K}}_1 \cup Z(G)$. Without loss of generality, we consider $Q := h_1 \bdot h_2 \bdot h_3$ with $h_1, h_2, h_3 \in {\mathcal{K}}_2 \cup {\mathcal{K}}_3 \cup {\mathcal{K}}_4$. Next, considering the previous subcase, we may assume that $U_1 \bdot Q$ does not contain a short $3$-sequence of thick type B, and $Q$ is not a short $3$-sequence. 

\smallskip

\noindent If $1 \in \Pi(Q)$, then we may assume that $h_1 h_2 \in Z(G)$, and then $U_1 \bdot h_1 \bdot h_2 \bdot g_{23} \bdot U_3 \bdot \dotsc \bdot U_9$ is an EGZ sequence, and we are done. So we assume that $1 \not\in \Pi(Q)$. Hence, if two of the terms of $Q$ are from the same $z$-class, they must be conjugate, or else they are from three distinct $z$-classes. In the first case we may assume that, either $h_1 h_2 \in {\mathcal{K}}_1$, or $h_1 h_2 h_3 \in {\mathcal{K}}_1$. Now if $h_1 h_2 \in {\mathcal{K}}_1$, we have $g_{11} g_{12} h_1 h_2 \in Z(G)$, since $U_1 \bdot Q$ does not contain a short $3$-sequence of thick type B. Then $g_{11} \bdot g_{12} \bdot h_1 \bdot h_2 \bdot g_{22} \bdot g_{23} \bdot U_3 \bdot \dotsc \bdot U_9$ is an EGZ sequence, and we are done. In case $h_1 h_2 h_3 \in {\mathcal{K}}_1$, one of the products $g_{11} h_1 h_2 h_3, g_{11} g_{12} h_1 h_2 h_3$ is central, and then $g_{11} \bdot h_1 \bdot h_2 \bdot h_3 \bdot g_{22} \bdot g_{23} \bdot U_3 \bdot \dotsc \bdot U_9$ (resp. $g_{11} \bdot g_{12} \bdot h_1 \bdot h_2 \bdot h_3 \bdot g_{23} \bdot U_3 \bdot \dotsc \bdot U_9$) is an EGZ sequence.  Now assume $h_1 \in {\mathcal{K}}_2, h_2 \in {\mathcal{K}}_3, h_3 \in {\mathcal{K}}_4$. In this case, we claim that the product two of the terms of $Q$ must belong to ${\mathcal{K}}_1$, and then we are done by arguments as given above. If this is not true, then 
\[
h_1 h_2 \equiv h_3, h_2 h_3 \equiv h_1, h_1 h_3 \equiv h_2 \hspace*{.2in} {\mathrm{mod}}~ Z(G),
\]
since $U_1 \bdot Q$ does not contain a short $3$-sequence of thick type B. But, then
\[
(h_1 h_2 h_3)^2 \equiv (h_1 h_2)(h_2 h_3)(h_1 h_3) \equiv h_1 h_2 h_3 \hspace*{.2in} {\mathrm{mod}}~ Z(G), 
\]
which implies that $h_1 h_2 h_3 \in Z(G)$, which contradicts that $Q$ is not a short $3$-sequence.

\bigskip

\noindent {\bf Subcase IC.} None of the $U_i$ are of thick type.

\bigskip

\noindent In this case, if $U_i$ contains non-central elements, it must be of thin type. If two of the $U_i$ are of thin type that contain terms from two distinct $z$-classes, then from Lemma \ref{twin-thick-type}(c), we have $1 \in \pi(U)$. If three of the $U_i$ are of thin type that contain terms from a single non-central $z$-class, letting these by $U_1, U_2, U_3$, we can permute the terms of these to form short $3$-sequences $W_1, W_2, W_3$ so that $W_1, W_2$ are of thick type A and $W_3$ is of central type, and the proof follows from subcase IB. Now suppose that $U_1, U_2$ are of thin type with terms from a non-central $z$-class ${\mathcal{K}}_1$, with $g_{13}, g_{23} \in Z(G)$, and the terms of $U \bdot (U_1 \bdot U_2)^{[-1]}$ are central. Now set $T := S \bdot U^{[-1]} = h_1 \bdot \dotsc \bdot h_6$. If two terms from $T$ are from ${\mathcal{K}}_1$, then we can replace $U_1, U_2$ by two short $3$-sequences with one of thick type A, and the other is either thick type A or of thin type. Now assume $h_1 \in {\mathcal{K}}_1$ and terms of $T \bdot (h_1)^{[-1]}$ are either central, or from the remaining non-central $z$-classes ${\mathcal{K}}_2, {\mathcal{K}}_3, {\mathcal{K}}_4$. Here we have $28$ elements from ${\mathcal{K}}_1 \cup Z(G)$ and using Lemma \ref{length-p3-subseqs}(iii), at least four elements of $T^{\prime} := T \bdot (h_1)^{[-1]}$ are from ${\mathcal{K}}_2 \cup {\mathcal{K}}_3 \cup {\mathcal{K}}_4$, and let $h_2, h_3, h_4, h_5$ be these elements. First, suppose $h_2 \bdot h_3 \bdot h_4 \bdot h_5$ contains a (minimal) central product subsequence $W$. If $\ell(W) = 2$, then we can replace $U_2$ by $W \bdot g_{23}$; the latter is a thin type sequence with terms from a $z$-class other than ${\mathcal{K}}_1$. If $\ell(W) = 3$, then $W$ is of thick type, and the argument follows from subcases IA and IB. If $\ell(W) = 4$, then $W$ must be an EGZ sequence and replacing $U_1 \bdot U_2$ by $W \bdot g_{13} \bdot g_{23}$ in $U$, the statement follows. Next, if $h_2 \bdot h_3 \bdot h_4 \bdot h_5$ does not contain a central product subsequence, then from Lemma \ref{Olson-minus1}, $h^{-1}_1 \equiv W^{\prime}$ modulo $Z(G)$ for some $W^{\prime} \mid h_2 \bdot h_3 \bdot h_4 \bdot h_5$ with $2 \leq \ell(W^{\prime}) \leq 4$. Now, making suitable replacements of terms of $U_1 \bdot U_2$ using $h_1 \bdot W^{\prime}$, as in the proof of Case IA, the statement would follow.  

\smallskip

\noindent Now we assume that $U_1$ is of thin type with non-central terms from ${\mathcal{K}}_1$ and $U_i ~(2 \leq i \leq 9)$ contain central terms. From Lemma \ref{length-p3-subseqs} (iii), it follows that the number of terms from the maximal subgroup ${\mathcal{K}}_1 \cup Z(G)$ that belongs to $T = S \bdot U^{[-1]}$ is at most $3$.  

\bigskip

\noindent First, assume that the three elements $z_1, z_2, z_3$ from $T$ are central and the remaining $h_1, h_2, h_3$ are non-central. Then the non-central elements are from ${\mathcal{K}}_2, {\mathcal{K}}_3, {\mathcal{K}}_4$. If these are from distinct $z$-classes, then from the previous case, we can assume that their product is non-central. But then, the product of two of these must belong to ${\mathcal{K}}_1$. If $h_1, h_2 \in {\mathcal{K}}_1$, then one of $g_{11} \bdot h_1 \bdot h_2$ and $g_{12} \bdot h_1 \bdot h_2$ is a short $3$-sequence of thick type and we are done by the previous case. Now suppose $h_1, h_2 \in {\mathcal{K}}_2$ and $h_3 \in {\mathcal{K}}_3$. Then $h_1 h_3, h_2 h_3 \in {\mathcal{K}}_4$. This implies that $h_1$ and $h_2$ are conjugate to each other. Then $h_1 h_2 h_3 \in {\mathcal{K}}_1$, and then one of the sequences $g_{11} \bdot h_1 \bdot h_2 \bdot h_3$ and $g_{12} \bdot h_1 \bdot h_2 \bdot h_3$ must have central product, and hence it is an EGZ sequence. This proves that all three elements of $T$ belong to the same $z$-class. Then, using Lemma \ref{twin-thick-type}(c), we are done. 

\smallskip 

\noindent Now, assume that $z_1, z_2$ from $T$ are central and the remaining elements $h_1, h_2, h_3, h_4$ are non-central. As argued  earlier, at most one of $h_i$ belongs to ${\mathcal{K}}_1$, and the remaining are from ${\mathcal{K}}_2, {\mathcal{K}}_3, {\mathcal{K}}_4$. The remaining arguments in this case are the same as in the previous paragraph. 

\smallskip

\noindent Finally, let at least five of the elements of $T$ be non-central. The proof of this part is similar to that in Subcase IA. 

\bigskip

\noindent {\bf Case II.} All terms of $U$ are central.

\bigskip

\noindent In this case, the sequence $S \bdot U^{[-1]}$ has at most one central element, and hence, there are either five or six terms of $S$ that are non-central. From Olson's Theorem \ref{Ols}, these non-central elements produce a subsequence $T^{\prime}$ of least length $l$ with $2 \leq l \leq 6$ and $\pi((T^{\prime})^{\ast}) \in Z(G)$. The proof now follows from the arguments as in Case I. \QED

\bigskip

\noindent {\bf Acknowledgement.} The authors would like to thank the anonymous referee for the suggestions for improvement of the paper and pointing out certain important errors in the previous version.

\bigskip

 
\bibliographystyle{plain} 
\bibliography{EGZ_Nonab_revised4_ArXiv}


\end{document}